\documentclass[12pt]{amsart}
\usepackage{mathrsfs, amssymb}
\usepackage[matrix,arrow,curve]{xy}

\makeatletter \@addtoreset{equation}{section} \makeatother

\newcommand{\FF}{\mathbb{F}}

\newcommand{\CC}{\mathbb{C}}
\newcommand{\QQ}{\mathbb{Q}}

\newcommand{\ZZ}{\mathbb{Z}}
\newcommand{\PP}{\mathbb{P}}

\newcommand{\OOO}{{\mathscr{O}}} 
\newcommand{\comp}{{{\,\scriptstyle{\circ}\,}}}
\newcommand{\down}[1]{\left\lfloor #1\right\rfloor}
\newcommand{\fr}[1]{\left\{ #1\right\}}

\newcommand{\muu}{{\boldsymbol{\mu}}}

\newcommand{\lin}{\text{---}}

\newcommand{\A}{\operatorname{A}}
\newcommand{\D}{\operatorname{D}}
\newcommand{\E}{\operatorname{E}}

\newcommand{\xref}[1]{{\rm \ref{#1}}}

\newcommand{\s}{\operatorname{s}}
\newcommand{\Pic}{\operatorname{Pic}}
\newcommand{\Sing}{\operatorname{Sing}}

\newcommand{\Bs}{\operatorname{Bs}}
\newcommand{\red}{\operatorname{red}}
\newcommand{\Supp}{\operatorname{Supp}}
\newcommand{\NE}{\overline{\operatorname{NE}}}
\newcommand{\Diff}{\operatorname{Diff}}
\newcommand{\T}{$\operatorname{T}$}

\renewcommand{\emptyset}{\varnothing}

\newtheorem{mtheorem}[equation]{}
\newtheorem{stheorem}{}[equation]
\newtheorem*{claim*}{Claim}

\theoremstyle{definition}
\newtheorem{mdefinition}[equation]{}
\newtheorem{sdefinition}[stheorem]{}

\title{A note on degenerations of del Pezzo surfaces}

\author{Yuri Prokhorov}
\thanks{
Partially supported by RScF grant no. 14-21-00053.
}

\address{
Steklov Institute of Mathematics, 8 Gub2kina street, Moscow 119991, Russia
\newline\indent
Department 
of Algebra, Faculty of Mathematics, Moscow State
University, Moscow 117234, Russia
\newline\indent
Laboratory of Algebraic Geometry, SU-HSE, 
7 Vavilova Str., Moscow 117312, Russia
}
\email{prokhoro@gmail.com}

\begin{document}

\maketitle
\section{Introduction}
\begin{mdefinition} 
This paper continues the classification of $\QQ$-Gorenstein degenerations of del Pezzo surfaces 
started in \cite{Manetti-1991}, \cite{Hacking-Prokhorov-2010}.
Let $\mathcal X\to \mathcal Z$ be a family of surfaces over a smooth curve $\mathcal Z$
such that a general fiber is a smooth del Pezzo surface and the special fiber 
$X:=\mathcal X_o$ is reduced, normal and has only quotient singularities.
Assume further that $\mathcal X$ is $\QQ$-Gorenstein and $-K_{\mathcal X}$
is ample over $\mathcal Z$. Such kind of families 
appear naturally in the three-dimensional minimal model program \cite{Kollar-ShB-1988}, \cite{Kollar-Mori-19988}
and in the study of certain moduli spaces  \cite{Hacking2004}, \cite{Hacking2012}.
It expected that the combinatorial structure of singularities of $X$
is related to exceptional vector bundles on smooth del Pezzo surfaces 
but this relation still looks mysterious (cf. \cite{Hacking2013}).

In this paper study  the special fiber $X$ of the above family 
under the condition that the Picard number $\rho(X)$ is large.
The case $\rho(X)=1$ was investigated completely in \cite{Hacking-Prokhorov-2010}.
Our main result is the following.
\end{mdefinition}

\begin{mtheorem}{\bf Theorem.} \label{main}
Let $X$ be a del Pezzo surface with only quotient singularities 
and let $\s(X)$ be the number of its non-Du Val points\footnote{For various definitions of Du Val (rational double) 
singularities we refer to 
\cite{Durfee1979}.}.
Assume that $X$ admits a $\QQ$-Gorenstein smoothing. Then 
\begin{enumerate}
\item 
$\s(X)\le \rho(X)+2$,
\item 
if $\s(X)=\rho(X)+2$, then $X$ is toric,
\item 
if $\s(X)=\rho(X)+1$, then $X$ admits an effective $\CC^*$-action.
\end{enumerate}
\end{mtheorem}

Similar to  \cite[Th. 1.3]{Hacking-Prokhorov-2010},
as a consequence of our techniques we verify a particular case 
of Reid's general elephant conjecture (cf.\cite[3.4B]{Reid-YPG1987}):

\begin{mtheorem}{\bf Theorem.}\label{thm-ge}
Let $f\colon \mathcal{X}\to \mathcal{Z}$ be a del Pezzo fibration over a smooth curve.
That is, $\mathcal{X}$ is a $3$-fold with terminal singularities, 
$f$ has connected fibers, and $-K_\mathcal{X}$ is ample over $\mathcal{Z}$.
Fix a point $o\in \mathcal{Z}$ and assume that the fiber $f^{-1}(o)$ is reduced, irreducible, normal, 
and has only quotient singularities.
Then, for some ample divisor $A$ on $\mathcal{Z}$, a general member $S\in |-K_\mathcal{X}+f^*A|$ 
is normal and has only Du Val singularities in a neighborhood of $f^{-1}(o)$.
\end{mtheorem}

Furthermore, we give a characterization of log surfaces that admit a $\CC^*$-actions
(Theorem \ref{theorem-toric})
and   establish the existence of $1$-complements 
on arbitrary del Pezzo surfaces with \T-singularities (Theorem \ref{th-1-compl}).
\par\medskip\noindent
\textbf{Acknowledgments.}

The work was written during my stay at the
International Centre for Theoretical Physics, Trieste.
I would like to thank ICTP for hospitality and support.

\section{Preliminaries}
\label{section-Preliminaries}
\begin{mdefinition}
Throughout this paper, we work over the field $\CC$ of complex numbers.
$\rho(X)$ denotes the Picard group of a variety $X$.
We use standard definitions, notation, and facts of the Minimal Model Program
\cite{Utah}, \cite{Kollar-Mori-19988}.

\begin{sdefinition}{\bf Proposition-Definition \cite[3.6]{Reid-YPG1987}, \cite[5.19]{Kollar-Mori-19988}.}
\label{def-canonical-index-one-covering}
Let $P\in X$ be the germ of a normal variety at a point $P$ and let $D$ be 
a Weil divisor on $X$. 
Assume  that $D$ is $\QQ$-Cartier at $P$, that is, $rD$ is
a Cartier divisor near $P$ for some positive integer $r$. 
Suppose that $r$ is the smallest such $r$. Then $r$ is called the \textit{index} of $D$.

There exists a covering $\pi: X^\sharp\to X$ which is Galois with 
group $\muu_r$ such that $ X^\sharp$ is normal, $\pi$ is etale over 
the the locus $X_0\subset X$ where $K_X$ is Cartier, and $P^\sharp:=\pi^{-1}(P)$ is a single point.
The divisor $K_{X^\sharp}=\pi^*K_X$ is Cartier.
Such a covering is called \textit{canonical index-one covering} of $P\in X$.
\end{sdefinition}
\end{mdefinition}

\begin{mdefinition}{\bf \T-singularities \cite{Kollar-ShB-1988}, \cite{Wahl-1981}.}
\label{Tdefclass}

\begin{sdefinition}{\bf Definition.}
Let $X$ be a normal surface such that $K_X$ is $\QQ$-Cartier.
We say that a deformation $\mathcal X/(0 \in S)$ over a germ $(0 \in S)$ with $0$-fiber $\mathcal X_0=X$ 
is \emph{$\QQ$-Gorenstein} if, locally analytically at every point $P \in X$,
$\mathcal X/S$ is induced by an equivariant deformation of 
the canonical index-one covering $(X^\sharp\ni P^\sharp)\to (X\ni P)$
(see \ref{def-canonical-index-one-covering}).
\end{sdefinition}

\begin{sdefinition}{\bf Definition \cite[Def.~3.7]{Kollar-ShB-1988}.} \label{defT}
Let $P \in X$ be a quotient singularity of dimension $2$.
We say $P \in X $ is a \emph{\T-singularity} if it admits a $\QQ$-Gorenstein smoothing.
That is, there exists a $\QQ$-Gorenstein deformation of $P \in X$ over a smooth curve germ
such that the general fiber is smooth.
\end{sdefinition}

\begin{stheorem}{\bf Proposition \cite[Prop.~3.10]{Kollar-ShB-1988}.} 
\label{Tclass}
A \T-singularity is either a Du~Val singularity or a cyclic quotient singularity of the form 
$\frac{1}{dn^2}(1,dna-1)$ for some positive integers $d$, $n$, $a$ with $\gcd(a,n)=1$.
\end{stheorem}
\end{mdefinition}

\begin{mdefinition}{\bf Noether's formula.}
For a \T-singularity $P \in X$, define
\[
\mu_P= 
\begin{cases}
r& \text{if $P \in X$ is a Du Val singularity of type $\A_r$, $\D_r$, or $\E_r$},\\
d-1&\text{if $P \in X$ is of type $\frac{1}{dn^2}(1,dna-1)$}. 
\end{cases}
\]
This number coincides with 
the \emph{Milnor number} of $P \in X$
\cite[Sec.~3]{Manetti-1991}.

\begin{stheorem}{\bf Proposition \cite{Hacking-Prokhorov-2010}.} \label{noether}
Let $X$ be a projective rational surface with \T-singularities. Then 
\[
K_X^2+\rho(X) + \sum_{P \in \Sing X} \mu_P = 10.
\]
\end{stheorem}

\begin{stheorem}{\bf Corollary.}\label{defspreservingrho}
Let $X$ be a del Pezzo surface with \T-singularities.
Then $K_X^2+\rho(X)\le 10$.
\end{stheorem}
\end{mdefinition}

\begin{mdefinition}{\bf Del Pezzo surfaces with \T-singularities.}
\begin{stheorem}{\bf Proposition \cite{Hacking-Prokhorov-2010}.} \label{dimension}
Let $X$ be a projective surface with \T-singularities such that $-K_X$ is nef and big. Then 
\[
\dim H^0(X,\OOO_X(-nK_X))= \frac12 n(n+1) K_X^2+1.
\]
\end{stheorem}
 \begin{stheorem}{\bf Corollary.} \label{dimension1}
Let $X$ be a projective surface with \T-singularities such that $-K_X$ is nef and big. Then 
\[
\dim |-K_X|= K_X^2>0.
\]
\end{stheorem}
\end{mdefinition}

\begin{mdefinition}{\bf Divisorial adjunction \cite[ch. 16]{Utah}.}
\label{def-diff}
Let $X$ be a normal variety and $S\subset X$ a reduced subscheme
of pure codimension one. Assume that the pair $(X,S)$ is
lc (log canonical \cite[2.10]{Utah}) in codimension two. Then there exists a naturally defined effective
$\QQ$-Weil divisor $\Diff_S(0)$, called the \emph{different},
such that
\[
(K_X+S)|_S=K_S+\Diff_S(0).
\]
Now let $B$ be a $\QQ$-divisor, which is $\QQ$-Cartier in
codimension two. Then the different for $K_X+S+B$ is defined by
the formula
\[
(K_X+S+B)|_S=K_S+\Diff_S(B). 
\]
In particular, if $B$ is a boundary and $(X,S+B)$ is lc in
codimension two, then $B$ is $\QQ$-Cartier in codimension two. Moreover, none of the
components of $\Diff_S(B)$ are contained in the singular locus of
$S$.
\end{mdefinition}

\begin{mdefinition}{\bf Classification of two-dimensional log canonical pairs
with reduced boundary \cite[ch. 3 \& Prop. 16.6]{Utah}, \cite[Th. 4.15]{Kollar-Mori-19988}.}
\label{2-dim-lc-pairs}
Let $P\in (X,C)$ be the germ of a two-dimensional log pair where 
$X$ is normal and $C$ is a (possibly reducible) reduced curve.
Assume that $(X,C)$ is lc. Then  one of the
following possibilities holds where all isomorphisms are 
isomorphisms of analytic germs:

\begin{sdefinition}\label{2-dim-lc-pairs-plt}
{\bf $(X,C)$  is plt {\rm (purely log terminal, \cite[2.13]{Utah}).}} Then 
\[
(X, C)\simeq(\CC^2, \{x_1=0 \})/\muu_{m}(1, a),\quad\text{with}\quad
\gcd(a, m)=1,\ m\ge 1,
\]
\[
\operatorname{Index}(K_X+C) =\operatorname{Index}(C) = m, \quad
\Diff_C(0)=(1-1/m)P.
\]
\end{sdefinition}

\begin{sdefinition}\label{2-dim-lc-pairs-lc-D}
{\bf $(X,C)$  is not plt and $C$ analytically reducible.} Then
\[
(X, C)\simeq(\CC^2, \{x_1x_2=0 \})/\muu_{m}(1, a),
\quad\text{with}\quad\gcd(a, m)=1,\ m\ge 1,
\]
\[
\operatorname{Index}(K_X+C) =1, \quad
\Diff_C(0)=0.
\]
\end{sdefinition}

\begin{sdefinition}\label{2-dim-lc-pairs-lc}
{\bf $(X,C)$  is not plt and $C$ analytically irreducible.} Then
\[
(X, C)\simeq(\CC^2, \{x_1x_2=0 \})/G,
\]
\[
\operatorname{Index}(K_X+C) =2, \quad
\Diff_C(0)=1,
\]
where $G\subset \operatorname{GL}_2(\CC)$ is a finite subgroup of dihedral type
without reflections (see \cite{Brieskorn-1967-1968} for the precise description
of $G$). 
\end{sdefinition}
\end{mdefinition}

\begin{mdefinition} {\bf $1$-complements.}
Let $X$ be a normal variety and let $D$ be a boundary on $X$
(an effective $\QQ$-divisor with coefficients $\le 1$).
Write $D=S+B$, where $S:=\down D$ (resp. $B:=\fr{D}$) is the integral 
(resp. fractional) part of $D$. A \emph{$1$-complement}
of $K_X+D$ is a divisor $D^+ \in |-K_X|$ such that $(X,D^+)$ is log canonical and 
$D^+ \ge S+\down{2B}$. In particular, if $D=0$, then a $1$-complement of 
$K_X$ is a divisor $D^+ \in |-K_X|$ such that $(X,D^+)$ is log canonical.
We say that 
the log divisor $K_X+D$ is \emph{$1$-complementary} if there exists a $1$-complement of $K_X+D$.

\begin{stheorem}{\bf Proposition  \cite[Prop. 4.3.2]{Prokhorov-2001}.}
\label{bir-prop}
Let $f\colon X\to Y$ be a birational contraction
and let $D$ be a boundary on $X$ such that
\begin{enumerate}
 \item 
\label{prop-i}
$K_X+D$ is nef over $Y$;
 \item 
\label{prop-ii}
the coefficients of $f_*D$ satisfy the
inequality $d_i\ge 1/2$.
\end{enumerate}
Assume that $K_Y+f_*D$ is $1$-complementary. Then so is $K_X+D$.
\end{stheorem}

\begin{stheorem}{\bf Proposition  \cite[Prop. 4.4.1]{Prokhorov-2001}.}
\label{prodolj}
Let $(X,D=S+B)$ be a log variety, where
$S:=\down{D}$ and $B:=\fr{D}$. Assume that
\begin{enumerate}
 \item 
$K_X+D$ is plt;
 \item 
$-(K_X+D)$ is nef and big;
 \item 
$S\ne 0$;
 \item 
\label{prop-2-iv}
the coefficients of $D=\sum d_iD_i$ satisfy the inequality $d_i\ge 1/2$.
\end{enumerate}
Further, assume that  there exists a
$1$-complement $K_S+\Diff_ S(B)^+$ of $K_S+\Diff_ S(B)$. Then 
 there exists a $1$-complement $K_X+S+B^+$ of
$K_X+S+B$ such that $\Diff_S(B)^+=\Diff_ S(B^+)$.
\end{stheorem}
\end{mdefinition}

\begin{mdefinition}{\bf Contractions of surfaces with  Du Val singularities.}
For convenience of the reader we state facts about MMP for Du Val surfaces.

\begin{sdefinition}{\bf Definition (see \cite{Morrison-1985}).}\footnote{For general definition see e.g. 
\cite[4.56]{Kollar-Mori-19988} or \cite[3.2]{Prokhorov-2001}.}
Let $y\in Y$ be a smooth point on a surface and let $(u,v)$ be a 
local coordinates near $y$.
A \textit{weighted blowup} with weights $(1,n)$ of a $y\in Y$
is the blowup $X\to Y$ of the ideal $(u,v^{n})$.
\end{sdefinition}

Clearly, a weighted blowup depends on $n$ and on the choice of coordinates.
For $n=1$ the above defined map is the usual blowup of $y\in Y$.
From easy local computations (see \cite[\S 1]{Morrison-1985}) we obtain the following.

\begin{stheorem}{\bf Lemma.}\label{Lemma-Morrison}
Let $y\in Y$ be a smooth surface germ and let $f: X\to Y$  be a weighted blowup with weights $(1,n)$,
$n\ge 2$.
Let $E=f^{-1}(y)$ be the exceptional divisor and let $\pi: \tilde X\to X$ be the minimal 
resolution.
Then the exceptional locus of the composition  $\tilde X\to Y$
is  a simple normal divisor whose dual graph looks as follows:
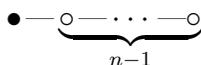
\[
\bullet \lin \underbrace{\circ \lin \cdots \lin\circ}_{n-1}
\]
where the vertex $\bullet$ corresponds to a $(-1)$-curve \textup(the proper transform of $E$\textup)
and the vertices $\circ$ correspond to $\pi$-exceptional $(-2)$-curves.
In particular, 
$X$ has exactly one singular point and this point is Du Val of type $\A_{n-1}$.
\end{stheorem}

\begin{stheorem}{\bf Corollary.}\label{Corollary-Morrison}
In the above notation we have $K_X\cdot E=-1$. 
\end{stheorem}

\begin{stheorem}{\bf Theorem (\cite[Theorem 1.4]{Morrison-1985}).}\label{theorem-Morrison}
Let $X$ be a surface with Du Val singularities and let $f: X\to Y$ be an
extremal Mori contraction. Let $E\subset X$ be the exceptional divisor and let 
$y:=f(E)$. Then $Y$ is smooth at $y$ and 
$f$ is a weighted blowup of $y\in Y$ with weights $(1,n)$ for some $n\ge 1$.
\end{stheorem}

\end{mdefinition}

\begin{mdefinition}{\bf Contractions of surfaces with \T-singularities.}
The following is the local variant of Theorem \ref{th-1-compl} below.
 \begin{stheorem}{\bf Proposition \cite[Prop. 4.7]{Prokhorov-2004b}.} 
\label{proposition-T-contractions}
Let $X$ be a surface with \T-singularities and let $f:X\to Y$ be a contraction such that
$-K_X$ is $f$-ample. Then, near each fiber $f^{-1}(y)$, $y\in Y$, there exists 
a $1$-complement of $K_X$. 
 \end{stheorem}

\begin{stheorem}{\bf Corollary.} 
\label{corollary-T-contractions}
Let $X$ be a surface with \T-singularities and let $f:X\to Y$ be a birational contraction such that
$-K_X$ is $f$-ample. If the fiber $f^{-1}(y)$ is not a point, then $y\in Y$ is a cyclic quotient
singularity \textup(or smooth\textup).
\end{stheorem}

\begin{proof}
Let $D\in |-K_X|$ be a $1$-complement of $K_X$ near the fiber. 
If $X$ has only Du Val singularities, then $Y$ is smooth at $y$ by Theorem \ref{theorem-Morrison}.
So we assume that $K_X$ is not Cartier near $f^{-1}(y)$. Then $D\neq 0$ and $D\cap f^{-1}(y)\neq \emptyset$.
Denote $D_Y:=f_*D$. 
Since $D\sim -K_X$ is $f$-ample, we have $\Supp(D)\not\subset  f^{-1}(y)$.
Hence $D_Y\neq 0$.
The pair $(Y,D_Y)$ is lc because so $(X,D)$ is.
Moreover, $K_Y+D_Y\sim f_*(K_X+D)\sim 0$.
By the classification \ref{2-dim-lc-pairs} the point $y\in Y$ is a cyclic quotient
singularity. 
\end{proof}

\begin{sdefinition} {\bf Warning.}
In general it is not true that the singularities of $Y$ are of type \T. 
\end{sdefinition}
\end{mdefinition}

\section{$\E$- and $\D$-singularities}
\label{section-singularities}

\begin{mtheorem}{\bf Proposition.}
\label{prop-DE}
Let $X$ be a del Pezzo surface with 
at worst quotient singularities and such that $\dim |-K_X|>0$.
Then $X$ has no Du Val points of type $\D_n$ or $\E_n$
contained in $\Bs |-K_X|$.
\end{mtheorem}

\begin{mtheorem}{\bf Corollary.}\label{corollary-DE}
Let $X$ be a del Pezzo surface with \T-singularities.
Then $X$ has no Du Val points of type $\D_n$ or $\E_n$
contained in $\Bs |-K_X|$.
\end{mtheorem}

\begin{mdefinition}
Let $P\in X$ be a Du Val point of type $\D_n$ or $\E_n$
such that $P\in \Bs |-K_X|$.
Let $D\in |-K_X|$ be a general member.
Write $D=\sum d_i D_i$, where the $D_i$ are prime divisors and $d_i>0$.
\end{mdefinition}

\begin{mtheorem}{\bf Lemma.}
Notation as above.
\begin{enumerate}
\item 
For any component $D_i$ of $D$ 
we have $D_i^2\ge 0$.
\item 
If $D_i^2= 0$ for some $D_i\subset D$, then 
$P$ is the only singular point of $X$ lying on $D_i$,
the pair $(X,D_i)$ is lc, and $d_i=1$.

\item 
All the components $D_i$ pass through $P$ and do not meet 
each other elsewhere.
\end{enumerate}
\end{mtheorem}
\begin{proof}
Let $D_i$ be a component passing through $P$. Assume that $D_i^2< 0$.
Then $D_i$ generates a birational extremal ray.
If $K_X$ is not Cartier along $D_i$, then by \cite[Cor. 4.3]{Prokhorov-2004b} $X$ has no Du Val 
points of type $\D_n$ or $\E_n$ on $D_i$, a contradiction.
If $K_X$ is  Cartier along $D_i$, then $X$ has only Du Val singularities in a neighborhood 
of $D_i$ and
we have a contradiction by Theorem  \ref{theorem-Morrison} (because our ray is $K_X$-negative).
This proves (i) modulo (iii).

Now assume that $D_i^2=0$ and $D_i\ni P$.
Then $D_i$ generates a contractible $K_X$-negative extremal face.
Thus there is a contraction $f:X\to Z$, where $Z$ is a smooth curve, such that $D_i=f^{-1}(z)_{\red}$
for some $z\in Z$. 
Since $D\in |-K_X|$ is a general member, 
the scheme fiber $f^*z$ is not contained in $D$. So, $f^*z\neq D_i$,
i.e. $f^*z$ is not reduced.
By \cite[Cor. 4.3]{Prokhorov-2004b}\ $K_X$ is Cartier along $D_i$.
Since $-K_X\cdot f^*z=2$, we have $-K_X\cdot D_i=1$ and $D_i$ is
a fiber of multiplicity $2$. Since $2D_i =f^*z$ is not contained in $D$, we have $d_i=1$.
Finally, by \cite[Prop. 7.1.3, Th. 7.1.12]{Prokhorov-2001} the pair $(X,D_i)$ is lc.
This proves (ii) modulo (iii).

Assume that (iii) does not hold. Then
there is a component $D_i$ such that 
$P\in D_i$ and $D_i\cap D_j\ni Q\neq P$
for some $j\neq i$ (because $\Supp D$ is connected). Put $D':=D-d_iD_i$.
By the above, $D_i^2\ge 0$. 
Hence the divisor $-(K_X+D_i+D')\sim (d_i-1)D_i$ is nef. Then by 
the adjunction we have
\[
\deg  \Diff_{D_i}(D')\le - \deg K_{D_i}\le 2.
\]
On the other hand,
$\Diff_{D_i}(D')\ge P+Q$. 
Hence, $\Diff_{D_i}(D')= P+Q$ and $(X,D_i+D')$ is lc near $D_i$
\cite[17.6]{Utah}. Since $P\in X$ is not a cyclic quotient
singularity,  the pair $(X,D_i)$ is strictly lc (i.e. lc but not plt) 
at $P$ (see \ref{2-dim-lc-pairs}).
In particular,
 no component of $D'$ pass through $P$.
Hence, $K_X+D_i+D'$ is not Cartier at $P$ (see \ref{2-dim-lc-pairs-lc-D}),
so $D\neq D_i+D'$ and $d_i>1$.
In this case, $D_i^2>0$ by (ii).
Then 
\[
(K_X+D_i+D')\cdot D_i=-(d_i-1)D_i^2<0
\]
and $\deg \Diff_{D_i}(D')<2$. The contradiction proves (iii). 
\end{proof}

\begin{proof}[Proof of Proposition {\rm \ref{prop-DE}}]
Let now $Q\in X$ be a non-Du Val point and let $D_i$ be a component
of $D$ passing through $Q$. Write $D=d_iD_i+D'$.

If  $D=D_i$, then
by \ref{2-dim-lc-pairs-plt} the pair $(X,D)$ is not plt at $Q$.
Hence, as above, $\Diff_{D}(0)= P+Q$.
By \cite[17.6]{Utah} $(X,D)$ is lc near $D$.
But then $K_X+D$ is not Cartier near $Q$, a contradiction.

Therefore, $D\neq D_i$.
Since $\dim |D|>0$, $D'\neq 0$ (and $D'$ has a reduced movable component). 
Note that the divisor $-(K_X+D_i+D')\sim (d_i-1)D_i$ is nef. Then by 
the adjunction we have
\[
\deg  \Diff_{D_i}(D')\le - \deg K_{D_i}\le 2.
\]
Since the coefficients of $\Diff_{D_i}(0)$ are $\ge 1/2$
and $(X,D_i)$ is not plt at $P$, $\Supp \Diff_{D_i}(D')=\{P,\, Q\}$.
Write $\Diff_{D_i}(0)=a_0P+bQ$ and $\Diff_{D_i}(D')=aP+bQ$.
Since $(X,D_i)$ is not plt at $P$, $a>a_0\ge 1$.
Hence $b<1$ and $(X,D_i)$ is plt at $Q$.
Thus, $b=1-1/m$ for some $m\ge 3$ (because 
$Q\in X$ is not Du Val of type $\A_1$)
and $a\le 4/3$.
If $(X,D_i)$ is not lc at $P$, then $a_0\ge 1+1/l$,
where $l$ is the minimal positive integer such that 
$l(K_X+D_i)$ is Cartier at $P$.
Recall that 
for the Weil divisor class group of a Du Val singularity $(X,P)$
we have 
\[
\begin{array}{c||ccccc}
(X,P)&
\D_{2n+1}&\D_{2n}&\E_6&\E_7&\E_8
\\[5pt]
\hline
\\[-5pt]
\operatorname{Cl}(X,P)&\ZZ/4\ZZ &\ZZ/2\ZZ\oplus \ZZ/2\ZZ& \ZZ/3\ZZ& \ZZ/2\ZZ&0
\end{array}
\]
(see, e.g., \cite{Brieskorn-1967-1968}).
So in our case we have $a_0\ge 5/4$ and 
$a\ge 5/4+1/4=3/2$, a contradiction. 

Thus we may assume that $(X,D_i)$ is  lc at $P$.
In particular, $(X,P)$ is of type $\D_n$.
Then $a_0=1$. Since $a\le 4/3$, we have only one possibility:
$b=2/3$, $a=5/4$, and $2D'$ is not Cartier at $P$. 
Moreover, $D'$ is irreducible (and reduced),
$(X,D')$ is not lc at $P$ (see \ref{2-dim-lc-pairs-lc-D}),
and so $\Diff_{D'}(0)\ge \frac54P$. 
Again by \ref{2-dim-lc-pairs-lc-D} $2D_i$ is Cartier at $P$.
Hence, 
\[
\Diff_{D'}(D_i)\ge \left(\frac54+\frac{d_i}2\right)P>2P,
\]
a contradiction.
\end{proof}

\section{Existence of $1$-complements}
\label{section-1-complements}

In this section we prove the following important fact (cf. \cite[Th. 7.1]{Hacking-Prokhorov-2010}).
\begin{mtheorem}{\bf Theorem.}
\label{th-1-compl}
Let $X$ be a del Pezzo surface with \T-singularities. 
Then there exists a $1$-complement of $K_X$. 
\end{mtheorem}

We need a few preliminary facts.

\begin{mdefinition}{\bf Definition \cite[\S 2]{Prokhorov-Shokurov-2009}.}
\label{CY->F}
Let $X$ be a normal projective variety. We say that
$X$ is \emph{FT} \textup(\emph{Fano type}\textup) if
there is a $\QQ$-boundary $\Delta$ such that 
$(X,\Delta)$ is a klt (Kawamata log terminal)  log Fano.
\end{mdefinition}

\begin{stheorem}{\bf Proposition \cite[\S 2]{Prokhorov-Shokurov-2009}.}
\label{CY-MMP}
Let $X$ be an FT variety. 
\begin{enumerate}
\item 
The Mori cone $\NE(X)$ is polyhedral and has contractible faces.
\item 
If $f: X\to Z$ be any contraction of normal varieties.
Then $Z$ is FT. In particular, the FT property is preserved under 
MMP.
\item
Let $\Xi$ be a boundary on $X$ such that $(X,\Xi)$ is lc and 
$-(K_X+\Xi)$ is nef.
Let $f: Y\to X$ be a birational extraction such that
$a(E,X,\Xi)< 0$ for every $f$-exceptional divisor $E$.
Then $Y$ is also FT.
\item 
Assume the LMMP in dimension $\dim X$.
Then the $D$-MMP works on $X$ with respect to any
divisor $D$.
\end{enumerate}
\end{stheorem}

\begin{mtheorem}{\bf Proposition.}
\label{prop-1-compl}
Let $(Y,C)$ be a log pair where $Y$ is an FT surface
and $C$ is an irreducible curve.
Assume that $(Y,C)$ is plt, $-(K_Y+C)$ is nef and big,
and $|-(K_Y+C)|\neq \emptyset$. Then one of the following holds:
\begin{enumerate}
\item 
$K_Y+C$ has a $1$-complement,
\item
$Y$ has three or four singular points on $C$ and either
\begin{enumerate}
\item 
$C^2<0$, $K_Y\cdot C\ge 0$, or
\item 
$\dim |-K_Y|=0$ and $-K_Y\sim bC$, $b\ge 2$.
\end{enumerate}
\end{enumerate}
\end{mtheorem}
\begin{proof}
First of all note that the curve $C$ is smooth 
(see \ref{2-dim-lc-pairs-plt}).
By Proposition \ref{prodolj} we can extend complements from 
$(C,\Diff_{C}(0))$ to $Y$.
Thus, for (i), it is sufficient to show existence of a $1$-complement of $K_{ C}+\Diff_{C}(0)$. 
Assume the converse and write
\[
\Diff_C(0)=\sum \left(1-\frac 1{m_i} \right)P_i, \quad
\deg K_C+\deg \Diff_C(0)=(K_X+C)\cdot C\le 0
\]
(see \ref{2-dim-lc-pairs-plt}). 
Thus, $\sum (1-1/m_i)\le 2$.
Since, by our assumption, the log divisor $K_{ C}+\Diff_{C}(0)$ is not $1$-complementary,
easy computations \cite[19.5]{Utah} 
show that $\Diff_{C}(0)$ is supported in 
three or four points $P_i$.
In particular, $\deg K_C<0$ and so $C\simeq \PP^1$.

Assume that $C^2\le 0$. Then $C$ generates an extremal face.
Since $Y$ is FT, this extremal face is contractible:
there is a contraction $\varphi: Y\to Y'$ such that $y:=\varphi(C)$ 
is a point. By Lemma \ref{lemma-2-points} below $K_Y\cdot C\ge 0$.
If $C^2=0$, then $Y'$ is a curve, $\varphi$ is a rational curve fibration, 
and $C=\varphi^{-1}(y)_{\red}$. 
In this case, $K_Y\cdot C<0$, a contradiction. 
Thus $C^2<0$ and we are 
in the case (iia).

Assume that $C^2>0$.
Let $D\in |-(K_Y+C)|$ be a general member. Write $D=aC+D'$,
where $a\ge 0$ and $C$ is not a component of $D'$.
If $D'=0$, we get case (iib). (Here $b=a+1\ge 2$ because
$K_X+C$ is not Cartier near singular points on $C$, see  \ref{2-dim-lc-pairs-plt}.) 
Thus we may assume that $D'\neq 0$.
Since the support of 
$$
(a+1)C+D'\in |-K_Y|
$$ 
is connected, $D'$ meets $C$. 
Further,
\begin{equation}
\label{eq-diff}
\deg \Diff_C( D')=-\deg K_C + (K_Y+ C+ D')\cdot C=2-(a-1)C^2\le 2. 
\end{equation}
By the above, $\Diff_C(D')$ has 
at least one point of multiplicity $\ge 1$
(and multiplicities of all points are $\ge 1/2$).
Since $\Diff_C(D')\ge \Diff_C(0)$, 
the only possibility is 
\[
\Diff_C(D')=P_1+\frac12 P_2+\frac12 P_3,
\]
where $P_1\in C\cap \Supp(D')$ and $P_2,\, P_3 \notin C\cap \Supp(D')$.
By \ref{2-dim-lc-pairs-plt} $K_Y+C+D'$ is not Cartier at $P_2$ and $P_3$.
Thus $K_Y+C+D'\not \sim 0$ and $a>1$.
On the other hand, $\deg \Diff_C(D')=2$, so by \eqref{eq-diff}
$a=1$, a contradiction.
\end{proof}

\begin{mtheorem}{\bf Lemma \cite[Prop. 7.1.12]{Prokhorov-2001}.}
\label{lemma-2-points}
Let $S\to Z$ be a $K$-negative extremal contraction
from a surface $S$ with log terminal singularities, 
where $Z$ is not a point.
Then $S$ has at most two singular points on each fiber.
\end{mtheorem}

\begin{proof}[Proof of Theorem {\rm \ref{th-1-compl}}]
Let $X$ be a del Pezzo surface with 
at worst quotient singularities and such that $\dim |-K_X|>0$
and let $D\in |-K_X|$ be a general member.
Take $t\in \QQ$ so that $(X,tD)$ is maximally lc.
If $t=1$, then $K_X+D$ is a $1$-complement.
So from now on we assume that $t<1$.

\begin{mdefinition}
Consider the case where $(X,tD)$ is plt. 
Write $tD=C+B$, where $C:=\down{tD}\neq \emptyset$
and $B$ is an effective fractional divisor.
Since $X$ is an FT variety,
we can run $-(K+C)$-MMP and obtain
\[
\varphi: X \longrightarrow \bar X.
\]
Since 
\[
-(K_X+C)\equiv B-(1-t)K_X, 
\]
all the contractions are $B$-negative.
Hence they are birational and we end up with a model 
$(\bar X,\bar C)$ such that $-(K_{\bar X}+\bar C)$ is nef.
We have 
\[
-(K_{\bar X}+\bar C)\equiv \bar B-(1-t)K_{\bar X}, 
\]
where
$-K_{\bar X}$ is ample and $\bar B:=\varphi_*B$ is effective.
Hence the divisor 
$-(K_{\bar X}+\bar C)$ is big. Further, 
\[
K_X+C+B\equiv -(1-t)D\equiv (1-t)K_X. 
\]
Hence all the contractions 
in $\varphi$ are $(K+C+B)$-negative.
Therefore, $(\bar X, \bar C+\bar B)$ is plt and so is 
$(\bar X, \bar C)$. 
So, $\bar B\neq 0$ and $\bar D:=\varphi_*D\neq \bar C$.
Apply Proposition \ref{prop-1-compl} to $(\bar X,\bar C)$.
The case (iia) does not occur because 
$-K_{\bar X}$ is ample and the case (iib) does not occur because
\[
\dim |-K_{\bar X}|\ge \dim |-K_X|>0.
\]
Hence, there exists a $1$-complement of $K_{\bar X}+\bar C$. 
By 
Proposition \ref{bir-prop} we can pull back $1$-complements 
from $\bar X$ to $X$.
\end{mdefinition}

\begin{mdefinition}
Now consider the case where $(X,tD)$ is not plt.
Put $B:=tD$.
Consider an inductive plt blowup \cite[Prop. 3.1.4]{Prokhorov-2001}\quad
$\delta: \hat X \to X$, that is, a birational extraction such that 
$\rho(\hat X/X)=1$ and 
\[
K_{\hat X}+\hat B+C=\delta^*(K_X+B),
\]
where $C$ is the (irreducible) exceptional divisor
and $\hat B$ is the strict transform of $B$.
Moreover, the pair $(\hat X, C +(1-\epsilon) \hat B)$ is plt 
for any $\epsilon >0$. 
Write 
\[
K_{\hat X}+\hat D+aC=\delta^*(K_X+D), 
\]
where $\hat D$ is 
the strict transform of $D$ and $a> 1$.
Then $\hat D+aC\in |-K_{\hat X}|$, so
$\dim |-K_{\hat X}|>0$. By Proposition \ref{CY-MMP} 
the variety $\hat X$ is FT. Run the $-(K+C)$-MMP.
As above all the contractions are $\hat B$-negative.
So we end up with a model $(\bar X,\bar C)$ where
$-(K_{\bar X}+\bar C)$ is nef and big (and $\bar C\neq 0$):
\[  
\xymatrix{
&\hat X\ar[dr]^{\varphi}\ar[dl]_{\delta}&
\\
X&&\bar X
}
\]
(the case $\hat X=\bar X$ is not excluded).
Since $\NE(\hat X)$ is polyhedral,
$-(K_{\hat X}+ C +(1-\epsilon) \hat B)$ is ample 
for some $0<\epsilon \ll 1$. Hence the plt property of 
the pair $(\hat X, C +(1-\epsilon) \hat B)$ is preserved.
In particular, $(\bar X,\bar C)$ is plt.
Apply Proposition \ref{prop-1-compl} to $(\bar X,\bar C)$.
The case (iib) does not occur because
\[
\dim |-K_{\bar X}|\ge \dim |-K_{\hat X}|>0.
\]
Assume that we are in the case (iia). 
Then $\bar C^2<0$ and $K_{\bar X}\cdot \bar C\ge 0$.
In particular, $\bar C$ is contractible: there is a contraction 
$\psi : \bar X\to \check X$ of $\bar C$, where $\check X$
is an FT surface. 
As in \cite[Proof of Th. 7.1]{Hacking-Prokhorov-2010} we see that
$\check P:=\psi(\bar C)\in \check X$ is a singular point and 
it is not a cyclic quotient singularity.
According to Zariski's main theorem the composition 
$\upsilon=\psi\comp\varphi\comp \delta^{-1}: X\dashrightarrow 
\check X$ is a morphism.
By Corollary \ref{corollary-DE} $\upsilon$ is not an isomorphism
(because $\delta(C)\in \Bs |-K_X|$).
Since $X$ is a del Pezzo, $\upsilon$ is 
a $K$-negative contraction.
On the other hand, by Corollary \ref{corollary-T-contractions}
$\check P$ is a cyclic quotient singularity,
a contradiction.
Therefore, the case (iia) does not occur 
and so there exists a $1$-complement of $K_{\bar X}+\bar C$. 
Now as above by 
Proposition \ref{bir-prop}  we can pull back $1$-complements 
from $\bar X$ to $X$.
\end{mdefinition}
\end{proof}

\section{Tori actions}
\label{section-Tori}
For a normal projective surface $X$ we denote by $\varrho(X)$ the \textit{numerical Picard number}, that is,
the rank of the group of Weil divisors modulo numerical equivalence.
Clearly, $\varrho(X)\ge \rho(X)$ 
and the equality holds if $X$ is $\QQ$-factorial.
For a $\QQ$-divisor $D=\sum d_iD_i$ on $X$ we denote
\begin{equation*}
\begin{array}{lll}
\| D \|&:=&\sum d_i,
\\[7pt]
\varsigma(X, D)&:= &\varrho(X)+2-\|D\|.
\end{array}
\end{equation*}
We say that a log pair $(X,D)$ is \textit{toric} if $X$ is a toric variety and 
$D$ is the (reduced) invariant boundary. 
We say that a log pair $(X,D)$ admits an \textit{effective $\CC^*$-action} if 
the variety $X$ admits such an action so that 
the divisor $D$ is $\CC^*$-invariant.

The statements (i) and (ii) of the following theorem
ware proved by Shokurov in much more general form \cite{Shokurov-2000}.
For the convenience of the reader we provide simplified complete proofs.

\begin{mtheorem}{\bf Theorem.}
\label{theorem-toric}
Let $(X, D)$ be a projective normal log surface 
such that $D$ is an integral \textup(effective\textup)
divisor, the pair $(X, D)$ is lc, and $K_X+D \sim 0$. Then
\begin{enumerate}
 \item 
$\varsigma(X, D)\ge 0$;
 \item 
if the equality holds, then $(X,D)$ is toric; 
 \item 
if $\varsigma(X, D)=1$, then $(X,D)$ admits an effective $\CC^*$-action. 
\end{enumerate}
\end{mtheorem}

\begin{sdefinition}{\bf Remark.}
Let $X$ be a projective normal surface and let $D\in |-K_X|$ be a divisor such that 
the pair $(X,D)$ is lc. Then
the property $\varsigma(X, D)= 0$ characterizes toric pairs.
On the other hand, the condition $\varsigma(X, D)\le  1$ is sufficient 
but not necessary 
for $(X,D)$ to admit an effective $\CC^*$-action.
For example, the product $C\times \PP^1$, where $C$ is an elliptic curve, 
admits an effective $\CC^*$-action
but for any $D\in |-K_X|$
we have $\varsigma(X, D)\ge 2$.
\end{sdefinition}

The rest of this section is devoted to the proof of Theorem
\ref{theorem-toric}.
We will use the following fact which 
is an easy consequence of the definition.

\begin{mtheorem}{\bf Lemma.}
\label{lemma-sigma-extremal}
Let $\varphi: Y'\to Y$ be a birational morphism of normal surfaces,
 and let
$D'$ be a reduced boundary on $Y'$. 
Denote by $N(\varphi, D')$ the number of $\varphi$-exceptional
curves that are not contained in the support of $D'$. 
Then
\[
\varsigma(Y',D')=
\varsigma(Y,g_*D')+N(\varphi, D').
\]
\end{mtheorem}

\begin{mdefinition}
Let $(X, D)$ be a projective log surface such that $(X, D)$ is lc, $K_X+D \sim 0$, and
\begin{equation}
\label{eq-main-contr}
\varsigma:=\varsigma(X, D)\le 1.
\end{equation}
Let $f: (X',D') \to (X,D)$ be a 
minimal dlt modification, that is, a birational map such that 
the log pair $(X',D')$ is dlt (divisorial log terminal \cite[2.13]{Utah}),
\[
K_{X'}+D'\sim f^* (K_X+D),\qquad f_*D'=D, 
\]
and any $f$-exceptional divisor 
has multiplicity $1$ in $D'$ (see e.g. \cite[Prop. 21.6.1]{Utah},
\cite[Prop. 3.1.2]{Prokhorov-2001}). By Lemma \ref{lemma-sigma-extremal}
\[
\varsigma(X', D')=\varsigma(X, D)=\varsigma\le 1.
\]
Hence, \eqref{eq-main-contr} holds for $(X',D')$.
Since $(X',D')$ is dlt, $X'$ is non-singular near $D'$ \cite[Prop. 16.6]{Utah}.
Moreover, $X'$ has at worst Du Val singularities outside of $D'$.

\begin{sdefinition}\label{MMP}
Run the $K$-MMP:
\[
(X',D')=(X^{(1)}, D^{(1)}) \overset{\varphi_1}\longrightarrow (X^{(2)}, D^{(2)}) \overset{\varphi_2} \longrightarrow \cdots 
\overset{\varphi_{l-1}}\longrightarrow (X^{(l)}, D^{(l)})= (Y, D_Y). 
\]
Let $E^{(i)}\subset X^{(i)}$ be the $\varphi_i$-exceptional
divisor. 
\end{sdefinition}
\end{mdefinition}

\begin{mtheorem}{\bf Claim.}\label{claim-Morrison}
For each $i=1,\dots,l$ we have 
\begin{enumerate}
\item 
$X^{(i)}$ has at worst Du Val singularities;
\item
$X^{(i)}$ is non-singular near $D^{(i)}$;
\item
$\varphi_{i}$ is the weighted blowup 
with weights $(1,n)$,  $n\ge 1$ of a smooth point $\varphi_{i}(E^{(i)}) \in X^{(i+1)}$;
\item 
$D^{(i)}$ is a simple normal crossing divisor;
\item
$E^{(i)}\cdot D^{(i)}=1$.
\end{enumerate}
\end{mtheorem}

\begin{proof}
One can prove (i)-(iii) by induction on $i$ using the following scheme:
\[
\text{ (i)${}_{i}$, (ii)${}_{i}$  $\Longrightarrow$  (iii)${}_{i}$, (i)${}_{i+1}$, (ii)${}_{i+1}$.}
\]
Indeed, if (i)${}_i$ holds, then $\varphi_i$ is a weighted blowup by Theorem \ref{theorem-Morrison}
and so $X^{(i+1)}$ is smooth at $\varphi_i(E^{(i)})$.

Since the pair $(X^{(i)}, D^{(i)})$ is lc, (ii) implies (iv) and
(v) follows from Corollary \ref{Corollary-Morrison}
because $ D^{(i)}\sim -K_{X^{(i)}}$.
\end{proof}

\begin{mtheorem}{\bf Claim.}
For each $i=1,\dots,l$ we have
$\varsigma(X^{(i)}, D^{(i)})\le \varsigma\le 1$,
in particular, $D^{(i)}\neq 0$.
\end{mtheorem}
\begin{proof}
It follows from Lemma \xref{lemma-sigma-extremal}. 
\end{proof}

By Lemma \xref{lemma-sigma-extremal} and because $X^{(i)}$ is non-singular near $D^{(i)}$,
on each step we have one of the following possibilities:
\begin{sdefinition}\label{claim-MMP-2-possibilities-i}
$\varsigma(X^{(i+1)}, D^{(i+1)})=\varsigma(X^{(i)}, D^{(i)})$,
$E^{(i)}\subset D^{(i)}$, 
and $\varphi_i$ is the usual blowup of a singular point of $D^{(i)}$;
\end{sdefinition}
\begin{sdefinition}\label{claim-MMP-2-possibilities-ii}
$\varsigma(X^{(i+1)}, D^{(i+1)})=\varsigma(X^{(i)}, D^{(i)})-1$, and
$E^{(i)}\not\subset D^{(i)}$.
\end{sdefinition}

\begin{stheorem}{\bf Corollary.}\label{corollary-MMP-2-possibilities}
Suppose that we are in the case \xref{claim-MMP-2-possibilities-i} above.
Furthermore suppose that $X^{(i+1)}$ admit an action of a 
connected algebraic group $G$ so that the boundary $D^{(i+1)}$
is $G$-invariant. Then the action lifts to $X^{(i)}$ so that $D^{(i)}$
is $G$-invariant.
\end{stheorem}

\begin{stheorem}{\bf Corollary.}\label{corollary-MMP-2-possibilities-a}
Suppose that we are in the case \xref{claim-MMP-2-possibilities-ii} above.
Furthermore suppose that $(X^{(i+1)},D^{(i+1)})$ 
is a toric surface.
Then the action of some one-dimensional 
subtorus $\mathcal T$ lifts to 
 $X^{(i)}$ so that $D^{(i)}$
is $\mathcal T$-invariant.
\end{stheorem}

\begin{proof}
Since $E^{(i)}\cdot D^{(i)}=1$ and $E^{(i)}\not\subset D^{(i)}$, the curve $E^{(i)}$ meets 
only one component $D^{(i)}_0\subset D^{(i)}$ so that $E^{(i)}\cdot D^{(i)}_0=1$.
Let $\pi: \bar X^{(i)}\to X^{(i)}$ be the minimal resolution near $E^{(i)}$.
By Lemma \ref{Lemma-Morrison}
the dual graph of $\bar X^{(i)}$ has the following form:
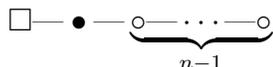
\[
\Box \lin \bullet \lin \underbrace{\circ \lin \cdots \lin\circ}_{n-1}
\]
where $\bullet$ corresponds to $E^{(i)}$, $\Box$ corresponds to $D^{(i)}_0$, 
and the vertices $\circ$ correspond to $\pi$-exceptional $(-2)$-curves.
Thus $\bar X^{(i)}$ is obtained from $X^{(i+1)}$ by making successive blowups of a fixed point
on the proper transform of $D^{(i)}_0$.
The stabilizer of this point is a one-dimensional subtorus in the big torus acting on $X^{(i+1)}$.
\end{proof}

\begin{mdefinition}
At the end of our MMP \ref{MMP} we get a log surface $(Y,D_Y)$
admitting a fiber type extremal $D_Y$-positive contraction $h: Y\to Z$. 
Moreover,
\begin{equation}
\label{eq-main-contr-Y}
\varsigma(Y, D_Y)\le \varsigma\le 1.
\end{equation}
Recall that $Y$ is non-singular near $D_Y$.
In particular, all the component of $D_Y$ are Cartier divisors.
Moreover, $K_Y+D_Y\sim 0$ and $Y$ has at worst Du Val singularities outside of $D_Y$.
\end{mdefinition}

\begin{mdefinition}\label{toric-del_Pezzo-case}
First we consider the case where $Z$ is a point.
Then $Y$ is a del Pezzo surface with at worst Du Val singularities and 
$\Pic (Y)\simeq \ZZ$. 
In particular, $\varrho(Y, D_Y)=1$.
Since $\|D_Y\|\ge 2$, the divisor $-K_Y$ is not a primitive element
of $\Pic(Y)$. In this case, $Y$ is either a projective plane $\PP^2$ or 
a singular quadric $\PP(1,1,2)$ (see e.g. \cite[Lemma 6]{Miyanishi-Zhang-1988}). Moreover
 we have one of the following:
\begin{sdefinition}
\label{end-3-components}
$Y\simeq\PP^2$, $D_Y=D_1+D_2+D_3$, where $D_i$ are lines in general position,
$\varsigma(Y, D_Y)=0$;
\end{sdefinition}
\begin{sdefinition}
\label{end-2-components-P2}
$Y\simeq\PP^2$, $D_Y=D_1+D_2$, where $D_1$ is a line and $D_2$ is a conic meeting 
$D_1$ transversely at two distinct points,
$\varsigma(Y, D_Y)=1$; 
\end{sdefinition}
\begin{sdefinition}
\label{end-3-components-P112}
$Y\simeq\PP(1,1,2)$, $D_Y=D_1+D_2$, where the class of $D_i$
generates $\Pic(Y)$ and again $D_1$ and $D_2$ meet each other transversely at two distinct points,
$\varsigma(Y, D_Y)=1$. 
\end{sdefinition}

\begin{claim*}
The pair $(Y, D_Y)$ is toric in the case \xref{end-3-components}
and admits an effective $\CC^*$-action in cases
\xref{end-2-components-P2} and \xref{end-3-components-P112}.
\end{claim*}

\begin{proof}
Modulo change of coordinates $x$, $y$, $z$  in $\PP^2$ or $\PP(1,1,2)$   we have 
\begin{itemize} 
\item[\ref{end-3-components}] $\Longrightarrow$ $D_Y=\{xyz=0\}$,
\item[\ref{end-2-components-P2}] $\Longrightarrow$ $D_Y=\{(xy-z^2)z=0\}$,
\item[\ref{end-3-components-P112}] $\Longrightarrow$  $D_Y=\{(xy-z)z=0\}$.
\end{itemize}
Then the statement of the claim is an easy exercise.
\end{proof}

In all cases \ref{end-3-components}--\ref{end-3-components-P112}
we have $\varsigma(Y, D_Y)\ge 0$.
This proves (i) of Theorem \ref{theorem-toric}. 
Moreover, if $\varsigma(X, D)= 0$,
then 
we are in the case \ref{end-3-components}.
In particular, $(Y,D_Y)$ is a toric surface.
Since $\varsigma(Y, D_Y)=\varsigma(X', D')=0$,
on each step of our MMP we have the
possibility \ref{claim-MMP-2-possibilities-i}.
By Corollary \ref{corollary-MMP-2-possibilities} both 
$(X',D')$ and $(X,D)$ are toric.

Now assume that $\varsigma(X, D)= 1$.
If moreover $\varsigma(Y, D_Y)=1$, then $(Y,D_Y)$ is 
of type \ref{end-2-components-P2} or \ref{end-3-components-P112}
and each step of our MMP is of type \ref{claim-MMP-2-possibilities-i}.
By Corollary \ref{corollary-MMP-2-possibilities} 
the action of the corresponding one-dimensional torus 
lifts to $X'$. Hence $(X',D')$ and $(X,D)$ admit effective $\CC^*$-actions.
Finally assume that $\varsigma(X, D)= 1$ and $\varsigma(Y, D_Y)=0$.
Then $(Y,D_Y)$ is toric and all but one steps of our MMP are of type 
\ref{claim-MMP-2-possibilities-i}.
As above we can apply Corollaries \ref{corollary-MMP-2-possibilities} 
and \ref{corollary-MMP-2-possibilities-a} to conclude that 
$(X,D)$ admits an effective $\CC^*$-action.
\end{mdefinition}

\begin{mdefinition}
Now consider the case where 
$Z$ is a curve. Then $Z$ is smooth and $h: Y\to Z$ is a rational curve fibration
with $\Pic(Y/Z)\simeq \ZZ$.
For a general fiber $F$ we have 
\[
D_Y\cdot F=-K_Y\cdot F=2. 
\]
Let $D_0$ be a $h$-horizontal component of $D_Y$. We claim that $D_0$ is a section.
Indeed, assume that $D_0$ is a double section. Then by the adjunction formula
\[
{D_0}\cdot (D_Y-D_0)=-{D_0}\cdot (K_Y+D_0)= -\deg K_{D_0}\le 2.
\]
Since $\|D_Y\|\ge 3$ and $D_Y$ is a simple normal crossing divisor, it 
has at least two 
vertical components $D_i$ with $D_0\cdot D_i=2$. Thus ${D_0}\cdot (D_Y-D_0)\ge 4$,
a contradiction. 

Hence, $D_0$ is a section. 
Then $D_Y$ has another $h$-vertical component $D_1$ which is also a 
section of $h$.
Since $D_0$ is a Cartier divisor, 
$h: Y\to Z$ is a smooth $\PP^1$-fibration.
If $D_0$ is not a rational curve, then as above by adjunction 
$\deg K_{D_1}= 0$ and $D_0$
is disjoint from $D_Y-D_0$. On the other hand, 
$D_Y-D_0$ has at least one $h$-vertical component, a contradiction. 
Hence, $D_0$ is a smooth rational curve and $Y$ is a Hirzebruch surface
$\FF_e$, $e\neq 1$. Let $\Sigma$ be the minimal section of $\FF_e$ and let 
$F$ be a fiber.
Since $D_i\cdot (D_Y-D_i)= -\deg K_{D_i}=2$ for each component 
$D_i\subset D_Y$, we have one of the following possibilities
(up to permutation of $D_0$ and $D_1$):

\begin{sdefinition}
\label{end-ruled-4-components}
$D_0\cdot D_1=0$, $\|D_Y\|=4$, $D_0=\Sigma$, $D_1\sim \Sigma+eF$,
$D_2$ and $D_3$ are distinct fibers, $\varsigma(Y, D_Y)=0$;
\end{sdefinition}

\begin{sdefinition}
\label{end-ruled-3-components}
$D_0\cdot D_1=1$, $\|D_Y\|=3$, $D_0=\Sigma$, $D_1\sim \Sigma+(e+1)F$,
$D_2$ is a fiber, $\varsigma(Y, D_Y)=1$.
\end{sdefinition}

Then we can complete the proof similar to \ref{toric-del_Pezzo-case}
by using the following.
\end{mdefinition}

\begin{claim*}
The pair $(Y, D_Y)$ is toric in the case \xref{end-ruled-4-components}
and admits an effective $\CC^*$-action in the case
\xref{end-ruled-3-components}.
\end{claim*}
\begin{proof}
The statement is obvious in the case $e=0$, so we assume that $e\ge 2$.
Let $\pi: (Y,D_Y)\to (Y',D_Y')$ be the contraction of the negative section.
Then $Y'$ is the weighted projective plane $\PP(1,1,e)$. We may assume that in suitable orbifold coordinates 
$x,\, y,\, z$ the boundary
$D_Y'$ is given by the equation 
$xyf_{e}(x,y,z)=0$ (resp. $xf_{e+1}(x,y,z)=0$) in the case \xref{end-ruled-4-components}
(resp. 
\xref{end-ruled-3-components}), where $f_{d}(x,y,z)$ denotes some polynomial of weighted degree $d$.

In the case \xref{end-ruled-4-components},
since $D_Y'$ is a simple normal crossing divisor outside of the origin $(0:0:1)$
and $D_1'$ does not pass through $(0:0:1)$,
the polynomial $f_{e}(x,y,z)$ contains $z$.
Then by a coordinate change we get $f_{e}(x,y,z)=z$. Hence $(Y',D_Y')$ is toric.
Since $\pi$ is the minimal resolution, the torus action lifts to $Y$.

Similarly, in the case \xref{end-ruled-3-components} 
$f_{e+1}(x,y,z)$ contains $zy$ (because $(Y',D_Y')$ is lc).
By a coordinate change we get $f_{e+1}(x,y,z)=zy+x^{e+1}$.
Then $(Y',D_Y')$ admits an $\CC^*$-action 
$(x,y,z)\longmapsto(x,\lambda y, \lambda^{-1}z)$.
\end{proof}

\section{Proof of main theorems}
\label{section-main-theorem}
Now Theorem \ref{main} is a consequence of the following.
\begin{mtheorem}{\bf Proposition.} \label{proposition-e-main}
Let $X$ be a projective normal surface
and let $\s(X)$ be the number of its points where $K_X$ is not Cartier.
Assume that $X$ has a $1$-complement $D\in |-K_X|$. Then 
\begin{enumerate}
\item 
$\s(X)\le \varrho(X)+2$,
\item 
if $\s(X)=\varrho(X)+2$, then $X$ is toric,
\item 
if $\s(X)=\varrho(X)+1$, then $X$ admits an effective $\CC^*$-action.
\end{enumerate}
\end{mtheorem}
\begin{proof}
By the classification of log canonical singularities of pairs \cite[Thm.~4.15]{Kollar-Mori-19988}, 
$D$ is a nodal curve,
and, at each singularity $P \in X$, either $D=0$ and $P \in X$ is a Gorenstein log canonical singularity, 
or the pair $P \in (X, D)$ is locally analytically isomorphic to the pair $(\frac{1}{n}(1,a),(uv=0))$ for some $n$ and $a$.
Moreover $D$ has arithmetic genus $1$ because $2p_a(D)-2=(K_X+D) \cdot D = 0$ (note that the adjunction formula holds
because $K_X+D$ is Cartier \cite[16.4.3]{Utah}). 
Thus $D$ is either a smooth elliptic curve, or a
rational curve with a node, or a cycle of smooth rational curves.

Let $\s'$ be the number of singular points of $X$ lying on $D$.
Then 
\[
\# \Sing(D) \ge \s'\ge \s(X).
\] 
By the above $\# \Sing(D)=\|D\|$. 
Then the assertion follows from Theorem \ref{theorem-toric}.
\end{proof}

The proof of Theorem \ref{thm-ge} is essentially the same as the proof of  \cite[Theorem 1.3]{Hacking-Prokhorov-2010}.

\section{Examples}
\label{section-Examples}
A natural way to produce examples of 
del Pezzo surfaces as in (iii) 
of Theorem \ref{main} is to use deformations:

\begin{mtheorem}{\bf Theorem \cite[Prop. 3.1]{Hacking-Prokhorov-2010}.} \label{unobs}
Let $X$ be a projective surface such that $X$ has only \T-singularities and 
$-K_X$ is nef and big. 
Then there are no local-to-global obstructions to deformations of $X$.
\end{mtheorem}

Thus we can start with some known examples and construct new ones by deforming 
their singularities. The behavior of the Picard number is described by 
Noether's formula \ref{noether} and by the following

\begin{mtheorem}{\bf Proposition \cite[Prop. 2.3]{Hacking-Prokhorov-2010}.} 
\label{Tdeformation}
Let $(P \in \mathcal X)/(0 \in S)$ be a $\QQ$-Gorenstein deformation 
of a \T-singularity $P\in X$ of type $\frac{1}{dn^2}(1,dna-1)$
and let $P_1,\dots,P_l$ be all the singular points of a fiber $\mathcal X_s$, $s\in S$.
Then the possible types of $P_1,\dots,P_l\in \mathcal X_s$  are as follows:
\begin{enumerate}
 \item[a)]
$A_{d_1-1},\ldots,A_{d_l-1}$ or 
 \item[b)] 
$\frac{1}{d_1n^2}(1,d_1na-1)$, 
$A_{d_2-1},\ldots, A_{d_l-1}$,
\end{enumerate}
where $d_1,\ldots,d_l$ is a partition of $d$. 
\end{mtheorem}

\begin{sdefinition}{\bf Remark.}
In the above situation, 
the case $\Sing(\mathcal X_s)=\emptyset$ is not excluded.
This is possible only if $d=1$ and 
in this case we put $l=1$.
\end{sdefinition}

\begin{stheorem}{\bf Corollary.}
Let $X$ be a projective surface with \T-singularities and 
let  $\mathcal X/(0 \in S)$ be a $\QQ$-Gorenstein deformation 
induced by a local deformation of one point $P\in X$. 
Then, in the notation of \xref{Tdeformation}, for a general fiber $\mathcal X_s$, $s\in S$ we have
\[
\rho(\mathcal X_s)-\rho(X)=l-1,
\]

\end{stheorem}

Now we can take  one of the toric surfaces with \T-singularities
and $\rho(X)=1$ 
described in \cite[\S 4]{Hacking-Prokhorov-2010} 
and deform it in a suitable way.
\begin{mdefinition}{\bf Example.}
 Take the weighted projective plane $X:=\PP(a^2,b^2,5c^2)$,
 where $a$, $b$, $c$,  
are given by the following Markov-type equation 
\[
a^2+b^2+5c^2=5abc
\]
(cf. \cite{Karpov1998}).
Then $X$ has three singular points which are of type \T{} and $K_X^2=5$.
More precisely, 
\[
\Sing(X)=\left\{\frac1{a^2} (b^2,5c^2),\ \frac1{b^2}(a^2,5c^2), \  \frac1{5c^2}(a^2,b^2) \right\} 
\] 
For the third point we have 
\[
 \frac1{5c^2}(a^2,b^2)=\frac1{5c^2}(1,5c\alpha-1),
\]
where $\alpha=ab\delta$ and $\delta$ is taken so that $a^2\delta\equiv 1\mod 5c^2$.
Thus deforming this point  to one of the following collection of singularities
\begin{itemize}
\item  
$\frac1{c^2}(1,c\alpha-1)$,\ $A_3$;
\item  
$\frac1{2c^2}(1,2c\alpha-1)$,\ $A_2$;
\item  
$\frac1{3c^2}(1,3c\alpha-1)$,\ $A_1$;
\item  
$\frac1{4c^2}(1,4c\alpha-1)$,
\end{itemize}
we get examples  of 
del Pezzo surfaces as in (iii) 
of Theorem \ref{main} with $K_X^2=5$, $\rho(X)=2$, $\s(X)=3$.
\end{mdefinition}

\def\cprime{$'$} \def\mathbb#1{\mathbf#1} \def\bblapr{April}

\end{document}